\documentclass[article,11pt]{amsart}    
\usepackage{amsmath}
\usepackage{amsfonts} 
\usepackage{amssymb}
\usepackage{dsfont}
\usepackage{paralist}
\usepackage{graphics} 
\usepackage{epsfig} 
\usepackage{graphicx}
\usepackage{epstopdf}

\usepackage[colorlinks=true]{hyperref}
\hypersetup{urlcolor=blue, citecolor=red}
\newtheorem{Proposition}{Proposition}[section]

\newtheorem{Lemme}{Lemma}[section]
\newtheorem{Theoreme}{Theorem}[section]

\newtheorem{Remarque}{Remark}

\def \R{\mathbb{R}}
\def \N{\mathbb{N}}

\def \finpv{\hfill $\blacksquare$}
\def \pv{{\bf{Proof.}}~} 

\def \ds{\displaystyle}

\renewcommand{\phi}{\varphi}
\newcommand\restr[2]{{
  \left.\kern-\nulldelimiterspace 
  #1 
  \vphantom{\big|} 
  \right|_{#2} 
}}

\begin{document}

\title[]{On the  blow-up for a Kuramoto-Velarde type equation}
\author{Oscar Jarr\'in}
\address{Escuela de Ciencias Físicas y Matemáticas, Universidad de Las Américas V\'ia a Nay\'on, C.P.170124, Quito, Ecuador}
\email{oscar.jarrin@udla.edu.ec}
\author{Gaston Vergara-Hermosilla}
\address{Laboratoire de Mathématiques et Modélisation d'Évry, CNRS UMR 8071,
		Université Paris-Saclay,  91025, \'Evry, France}
\email{(corresponding author) gaston.vergarahermosilla@univ-evry.fr}
\date{\today}


\maketitle

\begin{abstract}
It is known that  the  Kuramoto-Velarde equation is globally well-posed on Sobolev spaces in the case when the parameters $\gamma_1$ and $\gamma_2$ involved in the non-linear terms verify $ \gamma_1=\frac{\gamma_1}{2}$ or $\gamma_2=0$. In the complementary case of these parameters, 
the global existence or blow-up of solutions 
is a completely open (and hard) 
problem.
Motivated by this fact, 
in this work we consider
a non-local version of the  Kuramoto-Velarde equation.  This  equation allows us to apply a Fourier-based method  and, within  the framework $\gamma_2\neq \frac{\gamma_1}{2}$ and $\gamma_2\neq 0$, we show that large values of these parameters yield a blow-up in finite time of solutions in the Sobolev norm. 
\end{abstract}

\section{Introduction and main result}
The \emph{classical} Kuramoto-Velarde equation describes slow space-time variations of disturbances at interfaces, diffusion-reaction fronts
and plasma instability fronts \cite{Gar,Gar1}. The Kuramoto-Velarde   equation reads as:
\begin{equation}\label{Kuramoto-Velarde}
  \partial_t u
+
\partial_x^2 u
+
\partial_x^4 u
=
\gamma_1  (\partial_x u)^2
+
\gamma_2 \, u \partial^2_x u, 
\end{equation}
where $u:[0,+\infty)\times \R \to \R$ is a real-valued function and $\gamma_1, \gamma_2$ denotes constant parameters in $\R$. This equation also describes \emph{Benard–Marangoni cells}, which appear in the physical phenom of a  large surface tension on the interface  in a micro-gravity environment \cite{Coclite,Hym}. Precisely, the non-linearities $\gamma_1 (\partial_x u)^2$ and $\gamma_2 \, u \partial^2_x u$ model  pressure destabilization effects striving to rupture the interface \cite{Oer}. 

\medskip

A generalized version of the Kuramoto-Velarde equation (\ref{Kuramoto-Velarde}) is the \emph{dispersive} Kuramoto-Velarde equation, which, for a parameter $\alpha \geq 0$ reads as:
\begin{equation}\label{dispersive-Kuramoto-Velarde}
  \partial_t u
+
\partial_x^2 u
+
\alpha \partial^{3}_{x} u
+
\partial_x^4 u
=
\gamma_1  (\partial_x u)^2
+
\gamma_2 \, u \partial^2_x u. 
\end{equation}
This equation models  long waves on a viscous fluid flowing down an inclined plane  and
 drift waves in a plasma \cite{Cohen,Topper}. Mathematically speaking, the dispersive effects in this equation are given by the additional term $\alpha \partial^{3}_{x} u$. 

 \medskip

The analysis of  
 equation (\ref{dispersive-Kuramoto-Velarde}) posed in the whole line $\R$ was studied in \cite{Pilod}.
 In this equation, the effects of \emph{dissipative} terms  $\partial_x^2 u+\partial_x^4 u$ are (in some sense) stronger than the ones of the \emph{dispersive} term $\alpha\partial^{3}_{x} u$.
This fact allows the author of  \cite{Pilod}
to apply purely \emph{dissipative} methods in order to develop a well-posedness theory in Sobolev spaces. Precisely,   the \emph{local} well-posedness is proven in \cite[Theorem $1$]{Pilod}, in  the setting of non-homogeneous Sobolev spaces $H^s(\R)$, with $s>-1$. The value $s=-1$ is the optimal one in the sense that the flow map of  equation (\ref{dispersive-Kuramoto-Velarde}) is not a $\mathcal{C}^2$-function in $H^s(\R)$ with $s<-1$, see \cite[Theorem $3$]{Pilod}. 

\medskip

As usual, the existence time in the local well-posedness theory, which we shall denote by $T_0$, is driven by the size of the initial datum $u_0$:
\[ T_0 \lesssim \min\left( \frac{1}{(C\| u_0 \|_{H^s})^{1/\theta}}, 1 \right),\]
where $C$ and $\theta$ are positive quantities depending on $s>-1$. In addition, the solution $u$ also verifies the following regularity property $u \in \mathcal{C}((0,T_0], H^\infty(\R))$, where  we use the standard notation $\ds{H^\infty(\R)=\bigcap_{\sigma \geq s}H^\sigma (\R)}$.

\medskip

These same results hold for  equation (\ref{Kuramoto-Velarde}) (when $\alpha=0$) since, as mentioned above, the \emph{dispersive} term $\alpha\partial^{3}_{x} u$ does not play any substantial role in the well-posedness theory.  

\medskip

Always in the setting of the $H^s$-space (with $s>-1$), a remarkable feature of both equations (\ref{Kuramoto-Velarde}) and (\ref{dispersive-Kuramoto-Velarde})  is that the \emph{global} well-posedness issue \emph{strongly depends} on the  parameters $\gamma_1$ and $\gamma_2$ in the non-linear terms. In  fact, in \cite[Theorem $2$]{Pilod} it is proven that local solutions (obtained in \cite[Theorem $1$]{Pilod}) extend to global ones as long as:
\[ \gamma_2 = \frac{\gamma_1}{2} \quad \mbox{or} \quad \gamma_2=0.\]
In order to explain this fact, it is worth giving  a brief sketch   of the proof.
\begin{itemize}
    \item {\bf The case $\ds{\gamma_2 = \frac{\gamma_1}{2}}$}.  Global in time existence of solutions  essentially follows from a control on the $L^2-$norm of solutions. Remark that the second non-liner term can be rewritten as  $\ds{ u \partial^{2}_{x} u =  \frac{1}{2} \partial^{2}_{x}(u^2)-(\partial_x u)^2}$. Thereafter, since we have $u \in \mathcal{C}((0,T_0], H^\infty(\R))$ then $u$ solves (\ref{Kuramoto-Velarde}) and (\ref{dispersive-Kuramoto-Velarde})  in the classical sense, and taking the inner $L^2$-product of these equations with $u$ we get the following  energy estimate:
\begin{equation*}
\begin{split}
\frac{1}{2}\frac{d}{dt}\| u(t,\cdot)\|^{2}_{L^2}= &\, - \int_{\R}\partial^{2}_{x} u u dx - \alpha\int_{\R}\partial^{3}_{x} u u dx - \int_{\R}\partial^{4}_x u u dx\\
&\, +
 (\gamma_1-\gamma_2) \int_{\R} (\partial_x u)^2 u dx +\frac{\gamma_2}{2}\int_{\R} \partial^{2}_{x} (u^2) u dx.
 \end{split}
 \end{equation*}
Observe that it holds $\ds{\alpha\int_{\R}\partial^{3}_{x} u u dx=0}$. Moreover, by the Cauchy-Schwarz inequality and the relationship $ab-b^2 \leq \frac{1}{4}a^2$ (with $a,b\geq0$),  we have 
 \[ - \int_{\R}\partial^{2}_{x} u u dx  - \int_{\R}\partial^{4}_x u u dx \leq \frac{1}{4}\| u(t,\cdot)\|^{2}_{L^2}.\]
 On the other hand,  concerning the non-linear terms, by our assumption $\ds{\gamma_2=\frac{\gamma_1}{2}}$ and integrating by parts we obtain
 \begin{equation*}
  (\gamma_1-\gamma_2) \int_{\R} (\partial_x u)^2 u dx +\frac{\gamma_2}{2}\int_{\R} \partial^{2}_{x} (u^2) u dx = \left(\frac{\gamma_1}{2}-\gamma_2 \right)\int_{\R} (\partial_x u)^2 u dx =0.    
 \end{equation*}
Then, we are able to apply  Gr\"onwall inequality to get
 \begin{equation*}
   \| u(t)\|^{2}_{L^2} \leq \| u_0 \|^{2}_{L^2} e^{4t},  \end{equation*}
   which allows to extend the local solution to the whole interval of time $[0,+\infty)$.

\medskip

\item {\bf The case $\ds{\gamma_2=0}$}. 
Denoting by
$w=\partial_x u$, we begin by stressing the fact that this function  solves the equation
\begin{equation*}\label{Kuramoto-Velarde-3}
\partial_t w
+
\partial_x^2 w
+
\alpha\partial^3_x w
+
\partial_x^4 w
=
2 \gamma_1 w \partial_x w.
\end{equation*}
In addition, we obtain the standard non-linear transport-type term $w \partial_x w$, which verifies $\ds{\int_{\R} w \partial_x w\, w dx =0}$. Thus, the linear terms are estimated as above, and  we can obtain the control:
 \begin{equation*}
   \| w(t)\|^{2}_{L^2} \leq \| w_0 \|^{2}_{L^2} e^{ct}.  \end{equation*}
\end{itemize}

\medskip

To the best of our knowledge, the global existence of solutions in the complementary case of the parameters $\gamma_1$ and $\gamma_2$, {\it i.e.}
\begin{equation}\label{Framework}
\gamma_2 \neq \frac{\gamma_1}{2} \quad \mbox{and} \quad \gamma_2 \neq 0, 
\end{equation}   
remains  as a completely open (and hard) question. In this context, the main objective of this short article is to give some lights on the \emph{possible} blow-up of solutions in this case.

\medskip

The blow-up issue for non-linear PDEs is not trivial, and in some cases the original model must be slightly modified in order to  apply a rigorous method  to show blow-up in finite time of solutions. The main example is given in \cite{Montgomery-Smith}, where the author  studies the blow-up of solutions for a one-dimensional toy model of the three-dimensional Navier-Stokes equations. This idealistic model of the Navier-Stokes equations was generalized  in \cite{Gallagher} to the cases of two and three space dimensions, with the additional feature that the divergence-free condition is preserved.
See also  \cite{Cortez} and \cite{Chamorro}  for other examples concerning non-linear parabolic equations. 

\medskip

Inspired by this fact, we  introduce  the following non-local version of  equation (\ref{dispersive-Kuramoto-Velarde}):
\begin{equation}\label{Main-Equation}
\partial_t u
+
\partial_x^2 u
+
\alpha(-\partial^2_x)^{\frac{3}{2}} u
+
\partial_x^4 u
=
\gamma_1 \Big( (-\partial^{2}_{x} )^{\frac{1}{2}} u  
\Big)^2
+
\gamma_2 u \, \partial^{2}_x u.
\end{equation}
Here, for $\sigma \in \R$, we recall that the fractional Laplacian operator $(-\partial^{2}_{x})^{\frac{\sigma}{2}}$ (in dimension one) can be easily defined in the Fourier level by the symbol $|\xi|^\sigma$. Concerning the regularity,   equations (\ref{dispersive-Kuramoto-Velarde}) and (\ref{Main-Equation})  are quite similar since they have the same order in all the derivative terms. This fact yields that all the local well-posedness and regularity results mentioned above for the  equation (\ref{dispersive-Kuramoto-Velarde}) also hold for the equation (\ref{Main-Equation}). 

\medskip

In particular, for any initial data $u_0 \in H^s(\R)$, with $s>-1$, we can find a time $T_0>0$ (which depends on $\| u_0\|_{H^s}$ and $s$), and we can find a functional space  $\ds{X^{s}_{T_0}\subset \mathcal{C}([0,T_0], H^s(\R))}$, such that equation (\ref{Main-Equation}) has a unique solution $u \in X^{s}_{T_0}$. This functional space $X^{s}_{T_0}$ is defined as 
\begin{eqnarray*}
    X^{s}_{T_0}= \{ u \in  \mathcal{C}([0,T], H^s(\R)): \| u \|_{X^{s}_{T_0}}<+\infty\},
\end{eqnarray*}
with the norm 
\begin{equation*}
\| u \|_{X^{s}_{T_0}} = \sup_{0\leq t \leq T_0}\| u(t,\cdot)\|_{H^s}+\sup_{0<t\leq T_0} t^\frac{|s|}{4}\| u(t,\cdot)\|_{L^2}+\sup_{0<t\leq T_0} t^\frac{|s|+1}{4}\| \partial_x u(t,\cdot)\|_{L^2}.    
\end{equation*}
Here,  the second and the third expressions are designed to successfully control the non-linear terms in equation (\ref{Main-Equation}). The (unique) solution $u \in X^{s}_{T_0}$ is thus obtained by applying a standard fixed point argument to a \emph{equivalent} mild formulation (given in  expression (\ref{Mild-Equation}) below) of equation (\ref{Main-Equation}) . The proof of this  fact essentially follows the same arguments in the proof of \cite[Theorem $1$]{Pilod}, so we shall omit it.

\medskip

On the other hand, in contrast to equation (\ref{dispersive-Kuramoto-Velarde}), the non-local operators $(-\partial^2_x)^{\frac{3}{2}}$ and $(-\partial^{2}_{x} )^{\frac{1}{2}}$ in   equation (\ref{Main-Equation}) equation  have positive  symbols $|\xi|^3$ and $|\xi|$ in the Fourier variable. This fact  is one of our \emph{key-tools} to use a Fourier-based method in order to show blow-up in finite time of solutions to (\ref{Main-Equation}), when $\gamma_1$ and $\gamma_2$ verify the relationship (\ref{Framework}). 

\medskip

It is also worth mentioning this method  seems very difficult to be applied  in  equation (\ref{dispersive-Kuramoto-Velarde}):  the local operators $\partial^{3}_{x}$ and $\partial_x$  have  complex symbols $-i \xi^3$ and $i \xi$ in the Fourier variable, and consequently,  we loose all the sign information required to our method Fourier-based method. Nevertheless, equation (\ref{Main-Equation}) may be seen as a useful modification of equation   (\ref{dispersive-Kuramoto-Velarde}), in order to a better understanding on how the structure of non-linear terms:
\begin{equation*}
\gamma_1 \Big( (-\partial^{2}_{x} )^{\frac{1}{2}} u  
\Big)^2
+
\gamma_2 u \, \partial^{2}_x u,    
\end{equation*}
joint with sufficient conditions on the  parameters $\gamma_1$ and $\gamma_2$ (for instance (\ref{Framework})) can work together to yield the blow-up of solutions. 

\medskip

{\bf Main result.} We show that well-prepared initial data in equation (\ref{Main-Equation}) allows us to obtain an explicit blow-up time of solutions. For $s>-1$, we shall consider an initial datum $u_0 \in H^s(\R)$, explicitly  defined in the Fourier level by the expression:
\begin{equation}\label{Initial-data}
\widehat{u_0}(\xi) = \textcolor{black}{\eta}\, \mathds{1}_{
\left\{\left|\xi-3/2\right|<1 / 2\right\}}(\xi).
\end{equation}
Here, the quantity $\eta>0$ is suitable set to  verify the (technical) requirements:
\begin{equation}\label{Condition-eta}
\begin{cases}
    \eta^2  \gg 2^{1-2s}C_0, \quad \mbox{when} \ \ -1<s<\frac{1}{2}, \\
    \eta^2 \gg C_0, \quad \mbox{when} \ \ \frac{1}{2}\leq s,
\end{cases}  
\end{equation}
with a constant $C_0>0$ fixed big enough (see expression (\ref{C0}) below).  In this setting, our main result reads as follows:
\begin{Theoreme}\label{thm blow-up} Let $u_0$ be the initial datum defined  by (\ref{Initial-data}) and (\ref{Condition-eta}). Moreover, within the framework of the conditions (\ref{Framework}),  assume that $\gamma_1$ and $\gamma_2$ also verify 
\begin{equation}\label{Framework2}
\gamma_2 < 0 < \gamma_1,
\end{equation}
and 
\begin{equation}\label{Condition-gamma}
\frac{C_1}{2} \leq  \gamma_1 - \gamma_2,  \end{equation}
with $C_1>0$ big enough. 
Then, the  solution $u(t,x)$ of equation \eqref{Main-Equation} arising from $u_0$ blows-up at the time $\ds{T_*= \frac{2\ln(2)}{3}}$, and for $s>-1$ we have $\ds{\| u(T_*,\cdot)\|_{\dot{H}^s}=+\infty}$ . 
\end{Theoreme}

\medskip

The following comments are in order.  First, observe that our initial data $u_0$ defined in (\ref{Initial-data}) belongs to all the Sobolev spaces. We thus obtain the blow-up of solutions to equation (\ref{Main-Equation}) even in the case of a regular initial datum. 

\medskip

Conditions  (\ref{Framework2}) and (\ref{Condition-gamma}) on the parameters $\gamma_1$ and $\gamma_2$ are essentially technical requirements in order our method to work. However, they highlight that these parameters can strongly influence on the dynamics of solutions in equation (\ref{Main-Equation}).  

\medskip

Theorem \ref{thm blow-up} holds for the (non-local) Kuramoto-Velarde equation:
\begin{equation*}
 \partial_t u
+
\partial_x^2 u
+
\partial_x^4 u
=
\gamma_1 \Big( (-\partial^{2}_{x} )^{\frac{1}{2}} u  
\Big)^2
+
\gamma_2 u \, \partial^{2}_x u,   
\end{equation*}
which is a particular case of equation (\ref{Main-Equation}) when setting $\alpha=0$. Moreover, setting $\gamma_1=0$ and when the parameter $\gamma_2$ verifies $\gamma_2<0$ and $\frac{C_1}{2}\leq -\gamma_2$, this  theorem also holds for the original Kuramoto-Velarde equation (\ref{Kuramoto-Velarde}), with only the second non-linear term:
\begin{equation*}
  \partial_t u
+
\partial_x^2 u
+
\partial_x^4 u
=
\gamma_2 \, u \partial^2_x u.  
\end{equation*}
Here, it is interesting to observe  the strong effects of this non-linear term in equation (\ref{Kuramoto-Velarde}), which essentially block the global existence of solutions when the parameter $\gamma_2$ is negative and $|\gamma_2|$ is  large enough.

\section{Proof of Theorem \ref{thm blow-up}}

In this section we deal with the proof of Theorem \ref{thm blow-up}. For the sake of clearness, the proof of all our technical lemmas will be postpone to Appendix \ref{AppendixA}.\\

To begin, we consider  Duhamel's formula in order to recast equation (\ref{Main-Equation}) in the following (equivalent) form:
\begin{equation}\label{Mild-Equation}
u(t, x)
=
K(t, \cdot) 
* u_0(x)
+
\int_0^t K(t-\tau, \cdot) *\left(
\gamma_1 \Big( (-\partial^{2}_{x} )^{\frac{1}{2}} u  
\Big)^2 + \gamma_2 u \partial^2_x u \right)(\tau,\cdot) d \tau,
\end{equation} 
where the kernel $K(t,x)$ is explicitly  defined by the expression
\begin{equation}\label{Kernel}
K(t,x)
=
\mathcal{F}^{-1} 
\left(
e^{- 
\big( -\xi^2 + \alpha|\xi|^3 + \xi^4 \big) t}
\right)(x),
\end{equation}
with $\mathcal{F}^{-1}$ denoting the inverse Fourier transform in the spatial variable. Equation (\ref{Mild-Equation}) will be our \emph{key-element}  in the proof of Theorem \ref{thm blow-up}. In the sequel, we shall denote by $T_{max}>0$ the maximal time of existence of the solution $u$. 

\medskip

To continue, let us state the following fundamental lemma. 
\begin{Lemme}\label{Lema-Positivity}
  Let consider the initial data $u_0$ defined in (\ref{Initial-data}).  Moreover,  assume (\ref{Framework2}).   Then, the Fourier transform of the  solution $u$ to \eqref{Mild-Equation} is positive, {\it i.e.} $\widehat{u}(t,\xi)\geq 0$, for all $0\leq t \leq T_{max}$ and $\xi \in \R$. 
\end{Lemme}
As mentioned, this lemma is proven in  detail at Appendix (\ref{AppendixA}).  However, we highlight that condition (\ref{Framework2}) on $\gamma_1$ and $\gamma_2$ is one of the \emph{key-ingredients} in the proof. 

\medskip

With this result at hand, we are able to prove our main technical estimate:
\begin{Proposition}\label{Proposition-Lower-Bound} Define recursively  the following functions:
\begin{equation}\label{gk}
\begin{cases}\vspace{2mm}
\widehat{g_0}(\xi)=
\mathds{1}_{
\left\{\left|\xi-3/2\right|<1 / 2\right\}}(\xi), \\
\widehat{g}_n(\xi)=\widehat{g}_{n-1}(\xi) * \widehat{g}_{n-1}(\xi), \quad \mbox{for all $ n \in \mathbb{N}^*$},
\end{cases}
\end{equation}
and 
\begin{equation}\label{fk}
f_n(t)=e^{-\frac{3}{2}t(2^{n+4})} 2^{-5(2^n -1)}2^{5n}, \quad \mbox{for all $ n \in \mathbb{N}$}.
\end{equation}
Assume that $\gamma_1$ and $\gamma_2$ verify (\ref{Framework2}) and (\ref{Condition-gamma}). Moreover, let $T_*= \frac{2 \ln(2)}{3}$.  Then, the mild solution $u(t,x)$ of equation (\ref{Main-Equation}) can be bounded by below in the Fourier level as follows:
\begin{equation}\label{Lower-estimate}
   \widehat{u}(t,\xi)
    \geq \eta^{2^n} f_n(t)\,\widehat{g}_n(\xi),\quad  \mbox{for all $n\in \mathbb{N}$,  \  $t\geq T_*$,  and $\xi \in \R $}.
\end{equation}
\end{Proposition}

\noindent\pv  Remark that  mild solutions of equation \eqref{Main-Equation}  write  down in the Fourier variable as follows:
\begin{equation}\label{Mild-Fourier}
\begin{split}
   \widehat u(t, \xi)
= &\, 
e^{- t
\big( -\xi^2 + \alpha |\xi|^3 + \xi^4 \big) }  
 \widehat u_0(\xi)\\
 &\, +
\int_0^t 
e^{- t
\big( -\xi^2 + \alpha |\xi|^3 + \xi^4 \big) }  
     \left(
\gamma_1
(   |\xi| \widehat{u} * |\xi|\widehat{u}  )
-\gamma_2 \widehat{u}*(\xi^2 \widehat{u})
\right)(\tau,\xi)  
d \tau.
\end{split}
\end{equation}
Our starting point is to remark that assuming $0\leq \alpha \leq 1$, for all $\xi \in \R$ we have \begin{equation*}
-\xi^2 + \alpha |\xi|^3 + \xi^4 \leq -\xi^2 +  |\xi|^3 + \xi^4 \leq \frac{3}{2} \xi^4,
\end{equation*}
and then
\begin{equation*}
e^{- t
\big( -\xi^2 + \alpha |\xi|^3 + \xi^4 \big)} \geq e^{-\frac{3}{2}t \xi^4}.
\end{equation*}
\begin{Remarque}\label{Remark-alpha}
In the case $\alpha >1$ we  can  write
\[
-\xi^2 + \alpha |\xi|^3 + \xi^4 \leq  \frac{3}{2}(\alpha+1) \xi^4,\]
hence 
\[e^{- t
\big( -\xi^2 + \alpha |\xi|^3 + \xi^4 \big)} \geq e^{- \frac{3}{2}(\alpha+1) t \xi^4},\]
 and the subsequent estimates also depend on $\alpha$. Consequently, with a minor loss of generality, we can assume that $0\leq \alpha \leq 1$.
\end{Remarque}	

Moreover, recall that by (\ref{Initial-data}) we have $\widehat{u}_0\geq 0$, and by Lemma \ref{Lema-Positivity} and the relationship (\ref{Framework2})  we also have 
$\gamma_1
(   |\xi| \widehat{u} * |\xi|\widehat{u}  )
-\gamma_2 \widehat{u}*(\xi^2 \widehat{u}) \geq 0$.  With  these facts, we get back to (\ref{Mild-Fourier})  and we thus obtain
\begin{equation}\label{Mild-Fourier-2}
   \widehat u(t, \xi)
\geq   
e^{- \frac{3}{2} t \xi^4 
 }  
 \widehat u_0(\xi) +
\int_0^t 
e^{- \frac{3}{2} (t-\tau)\xi^4}  
     \left(
\gamma_1
(   |\xi| \widehat{u} * |\xi|\widehat{u}  )
-\gamma_2 \widehat{u}*(\xi^2 \widehat{u})
\right)(\tau,\xi)  
d \tau.
\end{equation}

We shall use this last expression to prove the estimate (\ref{Lower-estimate}). The proof follows from an induction process. 

\medskip

{\bf Step $n=0$}.  Here, we will prove the estimate: 
\begin{equation*}
   \widehat{u}(t,\xi) \geq  \eta f_0(t)\widehat{g}_0(\xi), \quad \mbox{for all $t\geq T_*$ and $\xi \in \R$}.  
\end{equation*}

In  expression (\ref{Mild-Fourier-2}), observe that the second term on right-hand side is positive. We thus write
\begin{equation*}
     \widehat u(t, \xi)
     \geq 
e^{- \frac{3}{2} t \xi^4}   
 \widehat u_0(\xi).
\end{equation*}
This fact, together with the definition of $\widehat u_0(\xi)$  and  $\widehat g_0(\xi)$ (given in expressions (\ref{Initial-data})  and   (\ref{gk}) respectively) yield the inequality 
\begin{equation*}
     \widehat u(t, \xi)
     \geq  \eta
e^{- \frac{3}{2}t \xi^4}        
\widehat{g}_0(\xi).
\end{equation*}

At this point, we state the following lemma: 
\begin{Lemme}\label{Lema-Coronas} Let $\widehat{g}_n(\xi)$ be defined as in (\ref{gk}).  For all 
 $n\in \N$ we have 
\begin{equation*}
\text{Supp}(\widehat{g}_n)
\subseteq \left\{\xi \in \mathbb{R}:   2^n<|\xi|< 2^{n+1}\right\}:= \mathcal{C}(2^n, 2^{n+1}).
\end{equation*}
\end{Lemme}

\medskip

Since $\text{Supp}(\widehat{g}_0)\subset \mathcal{C}(1,2)$, for all $\xi \in \R$ we get 
\begin{equation*}
 e^{- \frac{3}{2}t \xi^4}        
\widehat{g}_0(\xi) \geq   e^{- \frac{3}{2}t 2^4}        
\widehat{g}_0(\xi),
\end{equation*}
and thus, we can write
\begin{equation*}
  \widehat{u}(t,\xi) \geq \eta e^{-\frac{3}{2}t 2^{4}} \widehat{g}_{0}(\xi) = \eta  e^{-\frac{3}{2}t 2^{0+4}} \widehat{g}_{0}(\xi)= \eta  f_0(t) \widehat{g}_0(\xi).   
\end{equation*}

\medskip

{\bf Step $n \geq 1$}. Let assume that 
\begin{equation}\label{Assumption-n-1}
    \widehat{u}(t,\xi) \geq \eta^{2^{n-1}} f_{n-1}(t)\widehat{g}_{n-1}(\xi), \quad \mbox{for all $t \geq  T_*$ and $\xi \in \R$}.
\end{equation}
We shall prove that it holds for $n$. We get back to expression (\ref{Mild-Fourier-2}), hence we can write 
\begin{equation*}
\widehat u(t, \xi)
\geq    \int_0^t 
e^{- \frac{3}{2} (t-\tau)\xi^4}  
     \left(
\gamma_1
(   |\xi| \widehat{u} * |\xi|\widehat{u}  )
-\gamma_2 \widehat{u}*(\xi^2 \widehat{u})
\right)(\tau,\xi)  
d \tau.
\end{equation*}
By our assumption (\ref{Assumption-n-1}) we get
\begin{equation}\label{Estim01}
\widehat u(t, \xi)
\geq  \eta^{2^n}  \int_0^t 
e^{- \frac{3}{2} (t-\tau)\xi^4}  f^{2}_{n-1}(\tau) 
     \left(
\gamma_1
(   |\xi| \widehat{g}_{n-1} * |\xi|\widehat{g}_{n-1}  )
-\gamma_2 \widehat{g}_{n-1}*(\xi^2 \widehat{g}_{n-1})
\right)(\tau,\xi)  
d \tau.
\end{equation}
Thereafter, recall that by Lemma \ref{Lema-Coronas} we have $\text{Supp}(\widehat{g}_{n-1})\subset \mathcal{C}(2^{n-1},2^n)$. In particular,  the lower bound $2^{n-1}<|\xi|$  and the recursive definition (\ref{gk}) yield 
\begin{equation*}
|\xi| \widehat{g}_{n-1} * |\xi|\widehat{g}_{n-1}  \geq 2^{2(n-1)}  \widehat{g}_{n-1} * \widehat{g}_{n-1} = 2^{2(n-1)}  \widehat{g}_n, 
\end{equation*}
and 
\begin{equation*}
 \widehat{g}_{n-1}*(\xi^2 \widehat{g}_{n-1}) \geq    2^{2(n-1)}  \widehat{g}_n. 
\end{equation*}
Moreover, recalling that by (\ref{Framework2}) we have $0<\gamma_1$ and $\gamma_2<0$, in (\ref{Estim01}) we obtain
\begin{equation*}
\widehat u(t, \xi)
\geq  \eta^{2^n}  \int_0^t 
e^{- \frac{3}{2} (t-\tau)\xi^4}  f^{2}_{n-1}(\tau) (\gamma_1-\gamma_2) \, 2^{2(n-1)} \widehat{g}_n(\xi) d \tau. 
\end{equation*} 

At this point,  we should recall that hypothesis  (\ref{Condition-gamma}) we have $(\gamma_1 - \gamma_2) \geq C_1 2^{-1}$, thus, given a parameter $\beta >0$ (which we will be precise later) we write 
\[ (\gamma_1 - \gamma_2) \geq C_1 2^{-1} \geq C_1 2^{-1}  2^{-(n+1)\beta}. \]
Considering this in the last inequality, we get
\begin{equation*}
\begin{split}
    \widehat u(t, \xi)
\geq &\,  \eta^{2^n} C_1  \int_0^t 
e^{- \frac{3}{2} (t-\tau)\xi^4}  f^{2}_{n-1}(\tau) 2^{-1} 2^{-(n+1)\beta}\, 2^{2(n-1)} \widehat{g}_n(\xi) d \tau \\
=&\, \eta^{2^n}\, C_1\,  2^{-1} 2^{-(n+1)\beta+2(n-1)}\,  \int_0^t 
e^{- \frac{3}{2} (t-\tau)\xi^4}  f^{2}_{n-1}(\tau) \widehat{g}_n(\xi) d \tau.
\end{split}
\end{equation*}

\medskip

Now, recall that  by expression (\ref{fk}) for $0\leq \tau \leq t$ we have $f_{n-1}(\tau)\geq f_{n-1}(t)$. Then, we write

\begin{equation*}
\begin{split}
    \widehat u(t, \xi)
\geq &\, \eta^{2^n}\, C_1\,  2^{-1} 2^{-(n+1)\beta+2(n-1)}\, f^{2}_{n-1}(t)  \int_0^t 
e^{- \frac{3}{2} (t-\tau)\xi^4}   \widehat{g}_n(\xi) d \tau\\
=&\, \eta^{2^n}\, C_1\,  2^{-1} 2^{-(n+1)\beta+2(n-1)} \left( e^{-\frac{3}{2}  t (2^{(n-1)+4})} 2^{-5 (2^{n-1}-1)} 2^{5(n-1)} \right)^2 \int_0^t 
e^{- \frac{3}{2} (t-\tau)\xi^4}   \widehat{g}_n(\xi) d \tau\\
=&\, \eta^{2^n}\, C_1\,  2^{-1} 2^{-(n+1)\beta+2(n-1)} \left( e^{-\frac{3}{2} 2 t (2^{(n-1)+4})} 2^{-10 (2^{n-1}-1)} 2^{10 (n-1)} \right) \int_0^t 
e^{- \frac{3}{2} (t-\tau)\xi^4}   \widehat{g}_n(\xi) d \tau.
\end{split}
\end{equation*}
We must estimate the last integral above. By Lemma \ref{Lema-Coronas} we have the bound $|\xi|<2^{n+1}$ for all $\xi \in \text{Supp}(\widehat{g}_n)$. This fact and the inequality $t\geq  T_*$ (with $T_*= \frac{2\ln(2)}{3}$) yield
\begin{equation*}
  \begin{split}
\int_0^t 
e^{- \frac{3}{2} (t-\tau)\xi^4}   \widehat{g}_n(\xi) d \tau \geq &\, \int_0^t 
e^{- \frac{3}{2} (t-\tau)2^4(n+1)}   \widehat{g}_n(\xi) d \tau \\
=&\, \left(\frac{3}{2}\right)^{-1} 2^{-4(n+1)}  (1-e^{-\frac{3}{2} t 2^{4(n+1)}})  \widehat{g}_n(\xi)\\
\geq &\, \left(\frac{3}{2}\right)^{-1} 2^{-4(n+1)}  (1-e^{-\frac{3}{2}t})  \widehat{g}_n(\xi)\\
\geq &\, \left(\frac{3}{2}\right)^{-1} 2^{-4(n+1)} (1-e^{-\frac{3}{2} T_*}) \widehat{g}_n(\xi)\\
\geq &\,  \left(\frac{3}{2}\right)^{-1} 2^{-4(n+1)} 2^{-1} \widehat{g}_n(\xi).
  \end{split}  
\end{equation*}
Thus, in the previous estimate we obtain 
\begin{equation*}
\begin{split}
 \widehat u(t, \xi) \geq &\,   \eta^{2^n}C_1\,  2^{-1} 2^{-(n+1)\beta+2(n-1)} \left( e^{-\frac{3}{2} 2 t (2^{(n-1)+4})} 2^{-10 (2^{n-1}-1)} 2^{10 (n-1)} \right)  \left(\frac{3}{2}\right)^{-1} 2^{-4(n+1)} 2^{-1} \widehat{g}_n(\xi)\\
 =&\, \eta^{2^n} C_1 \left(\frac{3}{2}\right)^{-1}\,  2^{-1} 2^{-(n+1)(\beta+4)+2(n-1)} e^{-\frac{3}{2} 2 t (2^{(n-1)+4})}\left(  2^{-10 (2^{n-1}-1)} 2^{-1} \right)  2^{10 (n-1)}  \widehat{g}_n(\xi).
 \end{split}
\end{equation*}
Before continue, it is convenient to remark that 
\begin{equation*}
   2^{-(n+1)(\beta+4)+2(n-1)} = 2^{-(n+1)(\beta+4)+2(n+1-2)}= 2^{(n+1)(2-\beta-4)} 2^{-4}  ,
\end{equation*}
and
\begin{equation*}
 2^{-10 (2^{n-1}-1)} 2^{-1} = 2^{-5(2\times 2^{n-1}-2)}2^{-1} =2^{-5(2^n - 1 -1)} 2^{-1}=2^{-5(2^n -1)}2^4. 
\end{equation*}
Then, from the last inequality we can get 
\begin{equation*}
\begin{split}
 \widehat u(t, \xi) \geq &\, \eta^{2^n} C_1 \left(\frac{3}{2}\right)^{-1}\,  2^{-1} 2^{(n+1)(2-\beta-4)} e^{-\frac{3}{2} 2 t (2^{(n-1)+4})} 2^{-5(2^n -1)}2^{10 (n-1)}  \widehat{g}_n(\xi)\\
 =&\,   \eta^{2^n} C_1 \left(\frac{3}{2}\right)^{-1}\,  2^{-1 +2-\beta-4-10}  e^{-\frac{3}{2} 2 t (2^{(n-1)+4})} 2^{-5(2^n -1)}2^{n(2-\beta -4 +10)}  \widehat{g}_n(\xi)\\
 =&\,   \eta^{2^n} C_1 \left(\frac{3}{2}\right)^{-1}\,  2^{-13-\beta}  e^{-\frac{3}{2} 2 t (2^{(n-1)+4})} 2^{-5(2^n -1)}2^{n(8-\beta)}  \widehat{g}_n(\xi). 
 \end{split}
\end{equation*}
Here, we set $\beta=3$ and we obtain 
\begin{equation*}
\begin{split}
 \widehat u(t, \xi) \geq  &\,   \eta^{2^n} C_1 \left(\frac{3}{2}\right)^{-1}\,  2^{-16}  e^{-\frac{3}{2} 2 t (2^{(n-1)+4})} 2^{-5(2^n -1)}2^{5n}  \widehat{g}_n(\xi) \\
= &\, \eta^{2^n} C_1 \left(\frac{3}{2}\right)^{-1}\,  2^{-16}  e^{-\frac{3}{2}  t (2^{n+4})} 2^{-5(2^n -1)}2^{5n}  \widehat{g}_n(\xi)\\
=&\, \eta^{2^n} C_1 \left(\frac{3}{2}\right)^{-1}\,  2^{-16} f_n(t)  \widehat{g}_n(\xi). 
 \end{split}
\end{equation*}
Finally, assuming that the constant $C_1$ is large enough and it verifies $\ds{C_1 \left(\frac{3}{2}\right)^{-1} 2^{-16} \geq 1}$, we  obtain the wished estimate
\begin{equation*}
 \widehat u(t, \xi) \geq   \eta^{2^n} f_n(t)  \widehat{g}_n(\xi), \quad \mbox{for all $t\geq T_*$ and $\xi \in \R$}, 
\end{equation*}
and  Proposition \ref{Proposition-Lower-Bound} is proven. \finpv 
 
\medskip

The lower estimate (\ref{Lower-estimate}) obtained in Proposition \ref{Proposition-Lower-Bound} is now our key-tool to find a lower bound on the quantity $\| u(T_*,\cdot)\|_{\dot{H}^s}$, with $s>-1$.

\begin{Proposition}\label{Prop-lower-bound-norm} Let $s>-1$, and let $\gamma_1,\gamma_2$ numerical constants verifying (\ref{Framework2}) and (\ref{Condition-gamma}).  Then, there exists a generic constant $C>0$ such that the following estimate hold:
\begin{equation*}
\|  u(T_*,\cdot)\|^{2}_{\dot{H}^s} \geq
C\, \ds{\sum_{n=0}^{+\infty}} 2^{n(2s-1)} \left(
   \eta^2 e^{-\frac{3}{2}T_*(2^{5})} 2^{-10}\right)^{2^n}.
  \end{equation*} 
\end{Proposition}

\noindent\pv We write 
\begin{equation*}
\begin{split}
\| u(T_*,\cdot)\|^2_{\dot{H}^s}= &\, \int_{\R} |\xi|^{2s} |\widehat{u}(T_*,\xi)|^2 d \xi \geq  \sum_{n=1}^{+\infty} \int_{\mathcal{C}(2^n, 2^{n+1})} |\xi|^{2s} |\widehat{u}(T_*,\xi)|^2 d \xi \\
\geq &\, \sum_{n=1}^{+\infty} 2^{2sn} \int_{\mathcal{C}(2^n, 2^{n+1})}  |\widehat{u}(T_*,\xi)|^2 d \xi,
\end{split}
\end{equation*}
and  using the lower estimate (\ref{Lower-estimate}) we obtain
\begin{equation*}
 \sum_{n=1}^{+\infty} 2^{2sn}   \int_{\mathcal{C}(2^n, 2^{n+1})}  |\widehat{u}(T_*,\xi)|^2 d \xi \geq \sum_{n=1}^{+\infty}  2^{2sn} \,  \eta^{2^{n+1}} f^{2}_{n}(T_*)  \int_{\mathcal{C}(2^n, 2^{n+1})}\widehat{g}^2_{n}(\xi) d \xi.   
\end{equation*}
At this point, we shall need the following estimate:
\begin{Lemme}\label{Lema-Normas} Let $(\widehat{g}_n)_{n\in \mathbb{N}}$ be the family of functions defined in (\ref{gk}). Then, there exists a constant $C>0$ such that for all $n\in \mathbb{N}$, 
\[ \int_{\mathcal{C}(2^n, 2^{n+1})} \widehat{g}^2_{n}(\xi) d \xi \geq C 2^{-n}. \]
\end{Lemme}

\medskip

With this estimate,  and the definition of functions $f_n(T_*)$ (given in (\ref{fk})),  we get back to the previous inequality to get
\begin{equation*}
\begin{split}
  \sum_{n=1}^{+\infty}  2^{2sn} \eta^{2^{n+1}}  f^{2}_{n}(T_*)  \int_{\mathcal{C}(2^n, 2^{n+1})}\widehat{g}^2_{n}(\xi) d \xi 
 \geq &\,  C\, \sum_{n=1}^{+\infty}  2^{(2s-1)n} \eta^{2^{n+1}}  f^{2}_{n}(T_*)\\
 = &\, C\, \sum_{n=1}^{+\infty}  2^{(2s-1)n} \eta^{2^{n+1}}  \left( e^{-\frac{3}{2} 2 T_*(2^{n+4})}\,  2^{-10 (2^n -1)} 2^{10 n} \right)\\
 = &\, C\, \sum_{n=1}^{+\infty}  2^{(2s-1)n} \eta^{2^{n+1}}  \left( e^{-\frac{3}{2}  T_*(2^{n+5})}\,  2^{-10 \times 2^n} 2^{10} 2^{10 n} \right)\\
\geq &\, C\, \sum_{n=1}^{+\infty}  2^{(2s-1)n}  \left( \eta^2  e^{-\frac{3}{2}  T_*(2^{5})}\,  2^{-10} \right)^{2^n},  
 \end{split}
\end{equation*}
from where we obtain the desired inequality. Thus,
Proposition \ref{Prop-lower-bound-norm} is proven. \finpv 

\medskip

Now, we are able to finish the proof of Theorem \ref{thm blow-up}. For doing it, we need to consider the following cases of the parameter $s>-1$:

\medskip
\begin{itemize}
    \item {\bf Case  $-1<s<\frac{1}{2}$}.  In this case we have $2^{2s-1}<1$ and then $(2^{2s-1})^n \geq (2^{2s-1})^{2^n}$. Thus,  by Proposition \ref{Prop-lower-bound-norm} we get
    \begin{equation*}
    \|  u(T_*,\cdot)\|^{2}_{\dot{H}^s} \geq   C 
    \sum_{n=1}^{+\infty} \left(2^{2s-1}  
    \eta^2  e^{-\frac{3}{2}  T_*(2^{5})}\,  2^{-10}
    \right)^{2^n}.    
    \end{equation*}
Then, setting 
$\ds{\eta^2 \gg  2^{1-2s}C_0}$,  with  
\begin{equation}\label{C0}
 C_0= \left(  
    e^{-\frac{3}{2}  T_*(2^{5})}\,  2^{-10}\right)^{-1}, 
\end{equation}
and $T_*= \frac{2\ln(2)}{3}$, this series diverges and  we conclude the blow-up of the norm $ \| u(T_*,\cdot)\|_{\dot{H}^s}$.

\medskip

    \item {\bf  Case  $\frac{1}{2}\leq s$}. Remark that here we have $2^{n(2s-1)} \geq 1$. Considering  again Proposition \ref{Prop-lower-bound-norm}, the constant $C_0$ defined above, and setting $\ds{\eta^2 \gg  C_0  }$ we obtain
     \begin{equation*}
    \| u(T_*,\cdot)\|_{\dot{H}^s} \geq   C\, \sum_{n=0}^{+\infty}  \left( \eta^2  e^{-\frac{3}{2}  T_*(2^{5})}\,  2^{-10} \right)^{2^n},
    \end{equation*}
    where the last series is also divergent. 
\end{itemize}
With this we conclude the proof of Theorem \ref{thm blow-up}. \finpv 
 
\section{Appendix: Proof of technical lemmas}\label{AppendixA}

\subsection*{Proof of Lemma \ref{Lema-Positivity}} With  a little abuse of notation we write $u_0=K(t, \cdot) 
* u_0$, hence, we recall that  by considering $n\in \mathbb{N}$ and 
\begin{equation*}
    u_{n+1}(t,\cdot)
=
K(t, \cdot) 
* u_0(x)
+
\int_0^t K(t-\tau,\cdot) *\left( \textcolor{black}{
\gamma_1
\Big( (-\partial^{2}_{x} )^{\frac{1}{2}} u_n  
\Big)^2}+ \gamma_2 u_n \partial^2_x  u_n 
\right)
(\tau, \cdot) 
d \tau,
\end{equation*}
the Picard iteration scheme provides an unique solution  $u$ to the  problem 
\eqref{Mild-Equation},   
where 
\begin{equation*}
    u = \lim_{n\to +\infty} u_n \in X^s_{T_0}
    .
\end{equation*}
By taking the Fourier transform in the space variable of each iteration $u_n$ we obtain
{\color{black}
\begin{eqnarray*}
   \widehat u_{n+1}(t, \xi)
&   
=
&
   \widehat K(t, \xi) 
 \widehat u_0(\xi)
+
\int_0^t 
    \widehat K(t-\tau, \xi) 
  \left(
\gamma_1 
( |\xi|  \widehat{u}_n *  |\xi| \widehat{u}_n)-
\gamma_2 \widehat{u}_n * (\xi^2 \widehat {u}_n)
\right)(\tau,\xi) 
\ d \tau, 
\end{eqnarray*}
}
where 
$$
\widehat K(t, \xi) = 
e^{-t 
\big( -\xi^2 + \alpha |\xi|^3 + \xi^4 \big) }
.$$
By hypothesis we have $\widehat{u_0}(\xi)\geq 0$,  then the positivity of the right-hand side above carries on in the Picard iteration as long as (\ref{Framework2}) holds. Thus,   $u$ satisfies 
$\widehat{u}(t,\xi)\geq 0$,  for all  $t\geq 0$ and $\xi \in \R$.    \finpv

\subsection*{Proof of Lemma \ref{Lema-Coronas}} 
The proof follows from an induction
process.

\begin{itemize}
    \item {\bf Step $n=0$}. This case follows easily from the fact that, the support of $\widehat g_0$ is given by $$\left\{  \xi \in \R  : | \xi - \frac{3}{2} | < \frac{1}{2}  \right\} =  ] 1,2 [.$$
    
    \item {\bf Step $n \geq 1$}. Let us  start by assuming  $\text{Supp}(\widehat{g}_{n-1})
\subset \mathcal{C}(2^{n-1}, 2^{n}).$ Now, considering this fact and  
$$
\widehat{g}_{n}=\widehat{g}_{n-1}*\widehat{g}_{n-1} = 
\int_{\R} \widehat{g}_{n-1}(\rho)\widehat{g}_{n-1}(\xi-\rho)d\rho,
$$
we conclude  $\text{Supp}(  \widehat{g}_{n-1}(\xi - \rho) ) \subseteq \{ 2^{n-1} < \xi - \rho < 2^{n} \}$. On the other hand, from the induction hypothesis we also know that $\text{Supp}(  \widehat{g}_{n-1}(\rho) ) \subseteq \{ 2^{n-1} <  \rho < 2^{n} \}$.
Thus, gathering together these facts yield,
$$\text{Supp}( \widehat{g}_{n-1}(\rho)\widehat{g}_{n-1}(\xi - \rho) ) 
\subseteq \,\, ] 2^{n}, 2^{n+1}  [, 
$$
and then we conclude $\text{Supp}(\widehat{g}_n)
\subseteq  \mathcal{C}(2^n, 2^{n+1})$.

\end{itemize}
With this we conclude the proof. 
\finpv

\subsection*{Proof of Lemma \ref{Lema-Normas}} 
This fact can be deduced by induction.  In fact, for $k=0$ we directly  obtain $\ds{\left\|\widehat{g}_0\right\|_{L^1} = 1}$. Now, let us suppose that for $k\geq 1$ we have $\|\widehat{g}_{k-1}\|_{L^1} = 1$. Considering the definition of functions $\widehat{g}_{k}$ and the fact the they are positive, we obtain  
\begin{equation*}
\begin{aligned}
\left\|\widehat{g}_k\right\|_{L^1} 
&
=
\int_{\mathbb{R}}| \widehat{g}_{k-1}
\ast 
\widehat{g}_{k-1} |
(\xi)  d \xi
 \\
& =\int_{\mathbb{R}}\int_{\mathbb{R}} \widehat{g}_{k-1}(\rho) \widehat{g}_{k-1}(\xi-\rho) d \rho d \xi \\
& =\int_{\mathbb{R}} \widehat{g}_{k-1}(\rho) \left(\int_{\mathbb{R}}
\widehat{g}_{k-1}(\xi-\rho) d \xi \right) d \rho.
\end{aligned}
\end{equation*}
By noticing  that $\ds{\int_{\mathbb{R}} \widehat{g}_{k-1}(\xi-\rho) d \xi 
=
\left\|\widehat{g}_{k-1}\right\|_{L^1} 
=1}$, we deduce  $\ds{\left\|\widehat{g}_{k}\right\|_{L^1} =1}$. \finpv 

\vspace{1cm}

\paragraph{\bf Acknowledgements.} 
The authors warmly thanks Diego Chamorro for his helpful comments and advises.
The second author is supported by the ANID postdoctoral program BCH 2022 grant No. 74220003.

\medskip

\paragraph{{\bf Statements and Declaration}}
Data sharing does not apply to this article as no datasets were generated or analyzed during the current study.  In addition, the authors declare that they have no conflicts of interest, and all of them have equally contributed to this paper.

\end{document}